\documentclass[letterpaper, 11pt]{article}
\usepackage{amssymb,amsmath}
\usepackage{epsf}
\usepackage{times,bm,mathrsfs}  
\usepackage{amsmath}
\usepackage{amssymb}
\usepackage{amsfonts}
\usepackage{mathrsfs}
\usepackage{bbm}
\usepackage{color}
\usepackage{graphicx}
\newcommand{\dcal}{\mathcal{D}}

\newcommand{\pcal}{\mathcal{P}}

\newcommand{\real}{\mathbb{R}}
\newcommand{\intgr}{\mathbb{Z}}
\newcommand{\nintgr}{\mathbb{N}}

\newcommand{\eps}{\varepsilon}

\newcommand{\ind}{\mathbbm{1}}

\newcommand{\bdm}{\begin{displaymath}}
\newcommand{\edm}{\end{displaymath}}
\newcommand{\bea}{\begin{eqnarray*}}
\newcommand{\eea}{\end{eqnarray*}}
\newcommand{\bean}{\begin{eqnarray}}
\newcommand{\eean}{\end{eqnarray}}

\newcommand{\prob}{\mathbb{P}}
\newcommand{\expec}{\mathbb{E}}

\newcommand{\poly}{\mathrm{poly}}

\newcommand{\coop}{\mathrm{C}}
\newcommand{\defec}{\mathrm{D}}
\newcommand{\action}{\mathbf{A}}
\newcommand{\allcoop}{\mathbf{C}}
\newcommand{\alldefec}{\mathbf{D}}
\newcommand{\baction}{\mathbf{B}}
\newcommand{\vol}{\mathrm{vol}}
\newcommand{\bush}{\mathbb{S}}
\newcommand{\bfn}{\mathbf{N}}
\newcommand{\bfs}{\mathbf{s}}

\newtheorem{theorem}{Theorem}

\newtheorem{corollary}{Corollary}

\newtheorem{lemma}{Lemma}

\author{{\bf Elchanan Mossel \,\,\,\,\,\,\,\,\,\, S\'ebastien Roch} \\ 
Department of Statistics\\
University of California, Berkeley\\
Berkeley, CA 94720-3860\\
\texttt{\small \{mossel,sroch\}@stat.berkeley.edu}}
\title{Slow Emergence of Cooperation\\ for Win-Stay Lose-Shift
on Trees}

\begin{document}

\maketitle

\thispagestyle{empty}

\begin{abstract}
We consider a group of agents on a graph who 
repeatedly play the prisoner's dilemma game against their neighbors.
The players adapt their actions to the past
behavior of their opponents by applying the win-stay
lose-shift strategy. On a finite connected graph, it is easy to see that the
system learns to cooperate by converging 
to the all-cooperate state in a finite time. 
We analyze the rate of convergence in terms
of the size and structure of the graph.
[Dyer et al., 2002] showed that the system converges rapidly on the
cycle, but that it takes a time exponential in the size of the graph to converge 
to cooperation on
the complete graph. 
We show that the emergence of cooperation is exponentially slow in 
some expander graphs. More surprisingly, we show that 
it is also exponentially slow 
in bounded-degree trees, where many
other dynamics are known to converge rapidly. 
\end{abstract}

\bigskip


\clearpage

\setcounter{page}{1}

\section{Introduction}

We consider a group of agents arranged on the nodes of a graph who 
repeatedly play the prisoner's dilemma game against their immediate neighbors.
The players adapt their actions to the past
behavior of their opponents by applying the so-called win-stay
lose-shift strategy~\cite{NS93} which, as the name suggests, 
consists in changing strategy whenever the payoff is deemed unsatisfactory. 
This model has been studied in the artificial
intelligence literature~\cite{Ki95} as a simple example
of ``co-learning''~\cite{ST93, ST97}.
On a finite connected graph, it turns out
that the system converges 
to the all-cooperate state---the globally optimal state---in finite time. 
In this respect, this instance of the iterated prisoner's dilemma (IPD) game on a graph provides 
an interesting example of a system learning to behave optimally
by a mechanism that involves each agent applying independently a simple 
strategy---or rule of thumb---which takes into account only
the latest actions of its immediate neighbors. For related work,
see~\cite{FL98} and references therein. See also~\cite{Ax84} for
the evolutionary perspective.

In order to understand how persistent this ``emergence of cooperation'' phenomenon is, it is
crucial to analyze the rate of convergence to the all-cooperate
state. Where the convergence is rapid, one would expect to observe the
optimal, cooperation state in a practical system based on similar dynamics. 
On the other hand, where the convergence is slow,
one would rather expect that such a system would stagnate in a
suboptimal, metastable state where a nonnegligible fraction of agents defect. 
Rates of convergence for IPD were studied in~\cite{Ki95, DG+02}
where the structure of the graph was shown to be a determining factor.

In this paper, we show that IPD exhibits an exponentially slow
convergence to cooperation on expander graphs and bounded-degree
trees. Our result for bounded-degree trees is somewhat surprising. In
particular, it should be compared to the behavior of global reversible
dynamics on trees~\cite{BKMP05} where the convergence is always rapid.
Note however that this slow convergence on trees is not
unprecedented. Notably, the contact process, a
common model of infection, is slow to converge on trees when the infection rate
is large. See e.g.~\cite{Li99} and references therein. In fact,
our proof suggests that IPD behaves very much like the contact process.
Nevertheless, the analysis of non-reversible particle systems has been an open challenge in
the last two decades and we hope that the results obtained here can shed some more
light on how such systems can be tackled.

The proof of slow convergence we give here combines several ideas. The
main idea is to look at the process at the right space-time scaling.
This approach, commonly used in probability (e.g. in the analysis of interacting particle
systems~\cite{Li85}), allows us to analyze the rough behavior of IPD---defection
survives for long periods of time in zones
that are densely populated by defectors.
The main technical difficulty is to control the dependencies between
different regions and different times. Then the process is compared to a directed percolation process
(where the directed axis corresponds to the time axis in the
original process). Using contour arguments we show that the directed
percolation process survives for an exponential time. See~\cite{Du84}
for background on directed percolation.

\subsection{Definitions and Previous Work}

Recall that the prisoner's dilemma game (PD) is a bimatrix game
with the following payoff matrix for the row player (and similarly for the column player):
\begin{displaymath}
\left(\begin{array}{cc}
R & S\\
T & P 
\end{array}\right)
\end{displaymath}
where $T > R > P > S$ and $2R > T + S$. The first row (column) corresponds
to the \emph{cooperate} action and the second row (column) corresponds
to the \emph{defect} action. The global---or Pareto---optimum is for both
agents to cooperate.
However, for any given action of the column player, it is always in the row
player's advantage to defect (and similarly for the column player).

For an agent playing PD, a simple way to adapt
to her opponent's behavior is the so-called Win-Stay Lose-Shift strategy
(WSLS)~\cite{NS93}, also known as the Pavlov
strategy~\cite{Ki95,ST97}. This works as follows. Every time
the game is played,
if the agent's payoff is one of the two smaller payoffs, i.e. $P$ or $S$, then
she switches her action in anticipation for the next round of play,
otherwise she keeps the same action.

We now consider a repeated graphical version of PD which we will
refer to as IPD.
Let $G=(V,E)$ be a
finite graph with $n = |V|$. Each node, $v$, is an agent
to which we associate an action $A_t(v) \in \{\coop,\defec\}$
at time $t \in \real_+$. (As will become clear in later sections,
it is easier to consider the continuous-time
version of this problem.) Here $\coop$ stands for \emph{cooperation}
while $\defec$ stands for \emph{defection}. The initial state
is $A_0(v) = \defec$ for all $v \in V$. The agents repeatedly play
PD against their
immediate neighbors in the graph through the following
mechanism. Each edge $e \in E$ has an exponential clock,
i.e. we associate to each edge an independent Poisson
process $\{T_i(e)\}_{i\geq 1}$ where  
all inter-arrival times $T_{i+1}(e) - T_i(e)$ are
independent Exp(1) (with the convention $T_0 = 0$). 
Every time a clock rings, say at edge $e = (u,v)$,
the endpoint agents $u$ and $v$ play one round of PD using
their respective actions $A_t(u)$ and $A_t(v)$, assuming
the clock rings at time $t$. Then the two agents update
their state using WSLS. In other words, if a clock rings
on edge $e = (u,v)$ at time $t$, we witness the following
transition for $(A_t(u), A_t(v))$
\begin{eqnarray*}
(\coop,\coop) &\to& (\coop, \coop)\\
(\coop, \defec) &\to& (\defec, \defec)\\
(\defec, \coop) &\to& (\defec, \defec)\\
(\defec, \defec) &\to& (\coop,\coop).
\end{eqnarray*} 

This defines a stochastic process for the state of the system
$\action_t = (A_t(v))_{v\in V}$ with initial state
the all-defect state,
$\action_0 = \alldefec \equiv (\defec, \ldots, \defec)$. It is clear that,
given the above allowed transitions, the system has a unique fixed point, the
all-cooperate state $\allcoop \equiv (\coop,\ldots,\coop)$. In particular,
if $G$ is a finite connected graph with $n \geq 2$, we have
a.s. $\action_t \to \allcoop$ as $t \to +\infty$. The question
of interest is: how long does it take to reach $\allcoop$ on a given
graph. It was shown by~\cite{DG+02}---and previously conjectured in~\cite{Ki95}---that 
the time to the emergence of cooperation
depends crucially on the structure of the graph. 
Let $T_\allcoop$ be the stopping time at which $\action_t$
reaches $\allcoop$ for the first time. Below, \emph{with high probability}
(w.h.p.)
means with probability $1 - 1/\poly(n)$ where $\poly(n)$ 
increases polynomially with $n$.
In~\cite{DG+02},
the following two results are proved.
\begin{theorem}[\cite{DG+02}]
Let $G$ be a cycle on $n$ vertices. Then
w.h.p. $T_\allcoop = O(n \log n)$.
\end{theorem}
\begin{theorem}[\cite{DG+02}]\label{thm:complete}
Let $G$ be the complete graph on $n$ vertices. Then
w.h.p. $T_\allcoop = \Omega((1.1)^n)$.
\end{theorem}

\subsection{Our Results}

Given the previous theorems, it is natural to conjecture that
the time to the emergence of cooperation is governed by the connectivity
of the graph: a high connectivity, as in the complete graph, leads to slow convergence, 
while a low connectivity,
as in the cycle, leads to fast convergence.
Surprisingly, we refute this intuition with our main
result.
\begin{theorem}\label{thm:tree}
There is a constant $d$ so that for all
$n$ there is a $d$-regular tree with $n$ vertices
for which w.h.p. $T_\allcoop = \Omega(\rho^n)$
for some $\rho > 1$ depending only on
$d$.
\end{theorem}
To prove this result, we study IPD on ``linear trees.''
The main technical ingredient is a coupling with
oriented percolation. The proof of this theorem
is given in Section~\ref{sec:tree}.

Although the connectivity conjecture turns out to be
wrong in general, the following theorem, an extension of
the complete graph result of~\cite{DG+02},
shows that the intuition is partly correct in one direction.
Let $G$ be a graph with $n$ vertices. Let
$\alpha, \beta$ be two increasing functions of $n$ such that
for all $n$, $0 < \alpha(n) < \beta(n) < n$. 
Define the $(\alpha,\beta)$-expansion
constant $\rho_{\alpha,\beta}(G)$ of $G$ as 
\begin{equation*}
\rho_{\alpha,\beta}(G) = \min \left\{\frac{|E(U, U^c)|}{\vol(U)}\ :\ U \subseteq V,\, 
\alpha(n) \leq |U| \leq \beta(n)\right\},
\end{equation*} 
where $E(U,U^c)$ is the set of edges between $U$ and $U^c$, 
$\vol(U)$ is the sum of the degrees of the nodes in $U$, and
$|X|$ is the cardinality of $X$. 
\begin{theorem}\label{thm:expand}
Let $\eps > 0$. Let
$\alpha, \beta$ be two increasing functions of $n$ such that
for all $n$, $0 < \alpha(n) < \beta(n) < n$.
Let $G$ be a 
graph with $n$ vertices such that
$\rho_{\alpha,\beta}(G) > 1/2 + \eps$. Then
there is a constant $a > 1$ (depending only on $\eps$) such that
w.h.p. $T_\allcoop = \Omega(a^{\beta(n) - \alpha(n)})$ (for $n$ large enough).
In particular, if $\alpha,\, \beta$ are linear in $n$,
the emergence of cooperation is exponentially slow.
\end{theorem}
This follows from a martingale argument similar to 
that used in~\cite{DG+02} which is detailed
in Section~\ref{sec:expand}. 
Note that in Theorem~\ref{thm:expand}, 
in order to obtain slow convergence,
it suffices to have large expansion
for relatively small sets. 
In particular, the theorem applies to expander graphs such as 
random regular graphs~\cite{Ka95,FKS89}.

\section{Win-Stay Lose-Shift on Trees}\label{sec:tree}

\begin{figure}
\begin{center}
\input{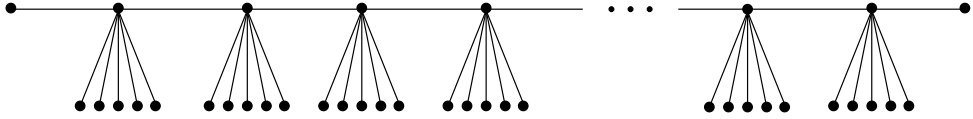}\caption{Caterpillar of degree $7$.}
\label{fig:cater}
\end{center}
\end{figure}

In this section, we analyze IPD on caterpillar trees
of degree $d$. 
We define an {\em (n,d)-caterpillar},
denoted $\bush_d^n$, 
to be a tree with the following property: 
the subtree induced by the internal nodes is a path
containing $n$ nodes all of which have degree $d$. 
See Figure~\ref{fig:cater}.
Our main
result, Theorem~\ref{thm:tree}, is that cooperation
is slow to emerge on caterpillars.
The proof of Theorem~\ref{thm:tree} follows from a series
of stochastic domination arguments. 
We now briefly outline the main steps of the proof.

\begin{enumerate}

\item \textbf{Star Dynamics via Biased Random Walk.} 
The first step is to analyze the behavior of a single star.  The
main point here is that it takes the star with $d$ leaves an exponential number
of steps (in $d$) to move from the all-defect state to the all-cooperate
state. This is proved by comparing the process to a biased random
walk. This comparison also shows that
a star can go from a few defectors to linearly many in $\poly(d)$ time
with constant probability, and
that a small linear fraction of defectors grows with high probability within $\poly(d)$
steps.  Moreover, these claims can be established even if one
allows two of the nodes of the stars to have arbitrary values.

\item \textbf{Space-Time Scaling.} 
We think of a star as defecting if at least $d/4$ of its
leaves defect.
Then, we consider triplets of adjacent stars and say that
a triplet is defecting if at least one of its extremal
stars is defecting. (We actually work with triplets of
stars rather than pairs to help control dependencies.)
We scale time by looking at the process
every $\poly(d)$ steps. The random walk argument of the previous point
allows to show that defecting stars have a high probability---at least
$(1-\exp(-\Omega(d)))$---of remaining defectors after
the $\poly(d)$ time window. Moreover, a defecting star has a $1/\poly(d)$ probability 
of ``infecting'' neighboring stars during that time.
By iterating these observations $\poly(d)$ times---yet another time scaling---we show 
that a defecting triplet has a probability
$1-\exp(-\Omega(d))$ of ``infecting'' a neighboring triplet. (Neighboring triplets
are actually intersecting.)

\item \textbf{Percolation.} We may now look at the space-time diagram of defecting triplets 
and show that
it dominates a directed percolation with probability
$1-\exp(\Omega(-d))$ for edges to be open. The time axis of the original
process corresponds to the direction of propagation in the percolation process.
Finally, a contour argument allows to conclude that this percolation
survives for a time which is exponential in $n$, thus proving that the
convergence time of IPD on the caterpillar is itself exponential in $n$.

\end{enumerate}

\subsection{Star}

Let $G = \bush_d^{n}$. This graph is made of 
$n$ copies of $\bush_d^1$ (i.e. stars of degree $d$). 
Let $G'$ be any star in $G$.
Denote the root $0$ and the leaves $1,2,\ldots,d$.
A crucial property of stars is that cooperation is slow to emerge on them.
This follows from our next result. 
We single out nodes $1$ and $2$, which are defined to be the two
nodes that $G'$ shares with its neighboring stars.
(In the case of extremal stars, we just pick an arbitrary
node in addition to the node shared with the next star.) 
We call $1$ and $2$ the external vertices. 
We use the following notation: $a\land b = \min\{a,b\}$.
\begin{lemma}[Dynamics on Stars]\label{lem:star}
Consider the IPD chain $\{\action_t\}_{t\geq 0}$ on $G=\bush_d^n$ with $d > 15$.
Let $G'$ be an arbitrary star in $G$ with nodes denoted $0,\ldots,d$
($0$ being the root, and $1$ and $2$ being the
external vertices).
Let $M'$ be a positive integer and $g_0,g_1,g_2$ be three
increasing functions
of $d$ with $g_2(d) = d/3 - 2$ and $g_0, g_1$ satisfying
$1 < g_0(d) < g_1(d) < g_2(d)$ for all $d$. 
Let the initial
configuration be as follows.
On $G'$, nodes $3$ through $d - g_1$ are $\coop$ and 
nodes $d - g_1 + 1$ through $d$ are $\defec$. 
On all other nodes, including the root and external vertices
of $G'$, the initial state
is arbitrary. Define
\begin{equation*}
N_\defec = |\{i\in \{3,\ldots,d\}\ :\ A(i) = \defec\}|.
\end{equation*}
Let $T_g$ be the
first time $N_\defec = g$. 
Let $\Delta_2 = g_2 - g_1$, 
$\Delta_1 = g_1 - g_0$,
$\rho = \sqrt{9/8}$,
and $\mu = g_0 M'$.
Then, we have
\begin{equation}\label{eq:star}
\prob[T_{g_2} \geq (T_{g_0}\land M')] \leq  
2^{-\Delta_1}
+ \rho^{-\sqrt{\mu}/2}(\sqrt{2})^{\Delta_2} + 2^{-\mu/2}.
\end{equation}
Moreover, this bound applies simultaneously on all stars
independently from each other (possibly with different
choices of $g$'s).
\end{lemma}
\noindent {\bf Proof:} 
For this argument, we restrict ourselves to
what happens on $G'$ and do not refer
to any event involving the rest of $G$. We call a leaf edge with leaf state $\defec$ a
$\defec$-edge, and similarly for $\coop$.
The behavior of $N_\defec$ depends on the state at the root of $G'$.
When $A(0) = \coop$, nothing happens until a $\defec$-edge
is picked at which time
$A(0)$ becomes $\defec$ itself.
On the other hand, when $A(0) = \defec$, either
a $\coop$-edge is chosen in which
case $N_\defec$ may go up by 1 (or stay the same
if $1$ or $2$ is picked),
or a $\defec$-edge is chosen in which
case $N_\defec$ may go down by 1 (or stay the same
if $1$ or $2$ is picked)
and $A(0)$ becomes $\coop$.
Ignore the updates where nothing changes, 
i.e. when an edge $(\coop,\coop)$ is chosen.
In any configuration satisfying
$N_\defec \geq g_0$, there are at least
$g_0$ edges whose updates change
the configuration. 
Let $Q$ the number of such updates in time $M'$.
Then it follows that $Q$ is larger than a Poisson with mean $\mu = g_0 M'$. 
From the moment generating
function of the Poisson distribution (see e.g.~\cite{Du96}), we have the following
\begin{equation*}
\prob[Q \leq \sqrt{\mu}] 
= \prob[e^{\mu - Q} \geq e^{\mu - \sqrt{\mu}}]
\leq \frac{\expec[e^{\mu - Q}]}{e^{\mu - \sqrt{\mu}}}
\leq \frac{e^{\mu} e^{\mu (e^{-1} - 1)}}{e^{\mu - \sqrt{\mu}}}
\leq 2^{-\mu/2}.
\end{equation*}
Assume the event $\{Q \geq \sqrt{\mu}\}$ holds.
Also, note that
at most one out of 2 steps have $A(0) = \coop$.
(Remember that we ignore $(\coop,\coop)$ updates.)
Ignore the times with $A(0) = \coop$ as well, what remains is
an asymmetric
random walk (or rather a birth-and-death chain)
which does at least $\sqrt{\mu}/2$ steps before time $M'$. 
To bound the probability that $N_\defec$ goes up or down,
we use the fact that the chain starts with $g_1$ $\defec$'s
and is stopped when it reaches either $g_0$ or $g_2$ $\defec$'s.
By assumption, the probability that $N_\defec$ goes up when
$A(0) = \defec$ is at least $(d - 2 - g_2)/d$. Consider
the walk $\{S_k\}_{k\geq 0}$ on $\nintgr$ started
at $S_0 = g_1$ which goes up with probability
$p = (d - 2 - g_2)/d = 2/3$ and goes down with probability
$1 - p = 1/3$.
Let $T'_g$ be the time at which $S_k$ reaches
$g$.
For convenience, we assume that the process $\{S_k\}_{k \geq 0}$
is defined on all of $\intgr$
(even though outside the interval $[g_0,g_2]$ the bounds
used are not valid).
Then,
\begin{equation*}
\prob[T_{g_2} \geq (T_{g_0}\land M')\,|\,Q \geq \sqrt{\mu}]
\leq
\prob[T'_{g_2} \geq (T'_{g_0}\land \sqrt{\mu}/2)]
\leq \prob[T'_{g_2} \geq T'_{g_0}] + \prob[T'_{g_2} \geq \sqrt{\mu}/2].
\end{equation*}
By standard martingale results (see e.g.~\cite{Du96}), we have
\begin{equation*}
\prob[T'_{g_2} \geq T'_{g_0}]
= \frac{\phi(\Delta_2) - \phi(0)}{\phi(\Delta_2) - \phi(-\Delta_1)},
\end{equation*}
where
\begin{equation*}
\phi(s) = \left(\frac{1 - p}{p}\right)^s = 2^{-s}.
\end{equation*}
So,
\begin{equation*}
\prob[T'_{g_2} \geq T'_{g_0}]
= \frac{1 - 2^{-\Delta_2}}{2^{\Delta_1} - 2^{-\Delta_2}}
\leq 2^{-\Delta_1}. 
\end{equation*}
We also have
\begin{equation*}
\expec[\rho^{T'_{g_2}}]
= \left(\frac{1 - \sqrt{1 - 4p(1 - p)\rho^2}}{2 (1-p)\rho}\right)^{\Delta_2}.
\end{equation*}
The choice $\rho = \sqrt{9/8}$ gives
\begin{equation*}
\expec[\rho^{T'_{g_2}}]
= \left(\sqrt{2}\right)^{\Delta_2}.
\end{equation*}
By Markov's inequality,
\begin{eqnarray*}
\prob[T'_{g_2} \geq \sqrt{\mu}/2]
=
\prob[\rho^{T'_{g_2}} \geq \rho^{\sqrt{\mu}/2}]
\leq \rho^{-\sqrt{\mu}/2}\left(\sqrt{2}\right)^{\Delta_2}.
\end{eqnarray*}
Finally, putting everything together, we get $(\ref{eq:star})$.

The independence of the bound at each star in $G$ comes from the fact
that we use only events involving leaf edges of $G'$.

\noindent $\blacksquare$

The following corollary corresponds to the case where a star 
has initially only a few $\defec$'s. The result below implies that
after $M' =\poly(d)$ steps the star has $O(d)$ $\defec$'s with positive
probability.
\begin{corollary}[Defection Spreads on Stars]\label{cor:birth}
In the setup of Lemma~\ref{lem:star}, let $g_0(d) = 2$,
$g_1(d) = 3$ and $g_2(d) = d/3 - 2$. Then, for 
$M' = \omega(d^2)$ and
$d$ (constant) large enough, we have
\begin{equation*}
\prob[T_{g_2} \geq (T_{g_0}\land M')] \leq \frac{2}{3}.
\end{equation*}
\end{corollary}

The following corollary implies that a star with $O(d)$ $\defec$'s
still has $O(d)$ $\defec$'s after $\poly(d)$ steps, 
with high probability.
\begin{corollary}[Defection Survives on Stars]\label{cor:survival}
In the setup of Lemma~\ref{lem:star}, let
$M' = +\infty$, $g_2(d) = d/3 - 2$, $g_1(d) = d/3 - 3$, 
and $g_0(d) = d/4 - 3$. Then,
\begin{equation*}
\prob[T_{g_2} \geq T_{g_0}] \leq  
2^{-d/12}.
\end{equation*}
\end{corollary}

The following corollary implies that a star with $d/4$ $\defec$'s
reaches $d/3$ $\defec$'s after $\poly(d)$ steps, 
with high probability.
\begin{corollary}[Defection Boosting on Stars]\label{cor:boosting}
Let $\tau$ be a positive integer, not depending on
$d$.
In the setup of Lemma~\ref{lem:star}, let
$g_2(d) = d/3 - 2$, $g_1(d) = d/4 - 2 - \tau$, 
and $g_0(d) = d/5 - 2 - \tau$.
Then, for $M' = \omega(d^2)$ and $d$ large enough, we have
\begin{equation*}
\prob[T_{g_2} \geq (T_{g_0}\land M')] \leq  
3\, 2^{-d/20} \leq 2^{-d/21}.
\end{equation*}
\end{corollary}

\subsection{Star Triplets}

\begin{figure}
\begin{center}
\input{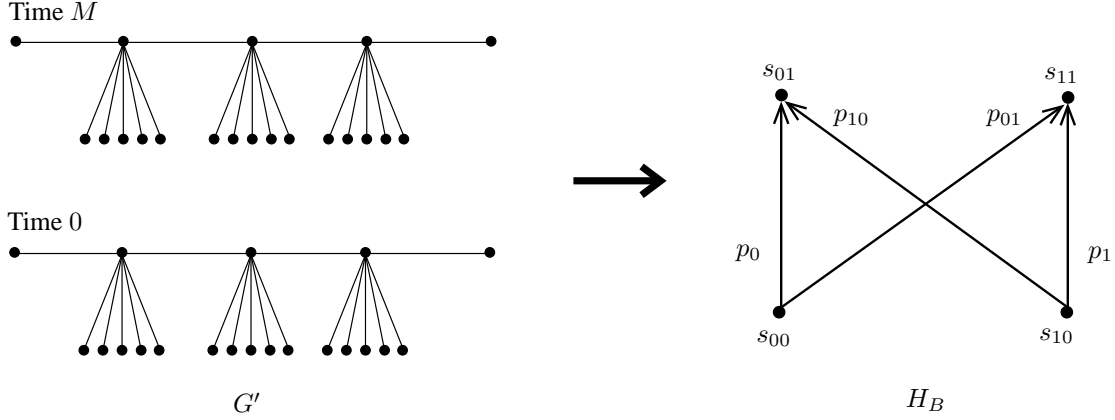}\caption{Reduction to percolation.}
\label{fig:perc1}
\end{center}
\end{figure}

The next step in the proof of Theorem~\ref{thm:tree} is to make
the connection between IPD and oriented percolation. Here
we show how a triplet of stars dominates the building
block of a percolation lattice. We use the following
oriented percolation. Consider four adjacent vertices of 
the regular lattice $\intgr^2$, say 
$v_{00} = (0,0)$, $v_{01} = (0,1)$, $v_{10} = (1,0)$ and 
$v_{11} = (1,1)$. Assume the nodes are connected
by four directed edges: $e_0 = (v_{00}, v_{01})$, 
$e_1 = (v_{10},v_{11})$, $e_{01} = (v_{00},v_{11})$, and
$e_{10} = (v_{10},v_{01})$. See Figure~\ref{fig:perc1}. 
Each edge is \emph{open} 
with respective probability $p_{0}$, $p_{1}$, $p_{01}$,
and $p_{10}$. The vertices have a state, denoted respectively
$s_{00}$, $s_{01}$, $s_{10}$, $s_{11}$, which
takes its value in $\{0, 1\}$. The
state $1$ ``travels''along the open edges, i.e.
if $e=(u,v)$ is an open edge and the state at $u$
is $1$ then the state at $v$ is also $1$. A vertex
is in state $1$ if and only if it is the terminal vertex of
an open edge with initial vertex in state $1$.
We denote
this four-node graph $H_B$. 

Now consider any triplet of adjacent stars inside $G = \bush_d^n$. 
Denote the stars $S_j$, $j=1,2,3$, with corresponding
edges $\{e_i^{(j)}\}_{j=1}^{d}$ and vertices $\{v^{(j)}_i\}_{i=0}^{d}$,
with the label $0$ corresponding to the root. We have the correspondence $e^{(1)}_2 = e^{(2)}_1$
and $e^{(2)}_2 = e^{(3)}_1$. 
We denote this subgraph---which is a copy of $\bush_d^3$---$G'$.
We are interested in the number of $\defec$'s on each 
star, excluding nodes $0$, $1$, and $2$ of each star, which we denote $\bfn_t = (N^{(1)}_t,
N^{(2)}_t, N^{(3)}_t)$.

The detailed behavior of $\bfn_t$ is rather intricate. We simplify the 
process by projecting it to a smaller space. Let
\begin{displaymath}
\sigma_d[N] 
= \left\{
\begin{array}{ll}
1, & \mathrm{if\ }N > d/4 - 2,\\
0, & \mathrm{if\ otherwise.}
\end{array}
\right.
\end{displaymath} 
Consider the random vector
\begin{equation*}
\tilde\bfs = (\tilde s_{00}, \tilde s_{01}, \tilde s_{10}, \tilde s_{11})
= \left(\sigma_d[N_0^{(1)}], \sigma_d[N_M^{(1)}], \sigma_d[N_0^{(3)}], \sigma_d[N_M^{(3)}]\right),
\end{equation*}
for some $M > 0$. The following lemma shows that for an appropriate
choice of $M$, $p_0$, $p_1$, $p_{01}$, and $p_{10}$, the vector $\tilde\bfs$
stochastically dominates
\begin{equation*}
\bfs = (s_{00}, s_{01}, s_{10}, s_{11}),
\end{equation*}
defined by the percolation above (with $s_{00} = \tilde s_{00}$ and
$s_{10} = \tilde s_{10}$).
\begin{lemma}[Connection to Percolation]\label{lem:triplet}
Consider the IPD chain $\{\action_t\}_{t\geq 0}$ on $G=\bush_d^n$ with $d > 15$.
Let $G'$ be an arbitrary triplet of adjacent stars in $G$.
Let $M = d^6$, $p_0 = p_1 = 1 - 2^{-d/30}$, and $p_{01} = p_{10} = d^{-10}$.
Then, for any initial configuration and $s_{00}, s_{10}$
such that $s_{00} = \tilde s_{00}$ and
$s_{10} = \tilde s_{10}$, we have that
$(\tilde s_{01}, \tilde s_{11})$ stochastically
dominates $(s_{01}, s_{11})$ for $d$ (constant) large enough. Moreover, the domination holds
for any number of (edge-)nonintersecting triplets simultaneously independently
from each other.
\end{lemma}
\noindent {\bf Proof:} 
The argument ignores any event outside $G'$. We consider three cases. 

\noindent 1) \textbf{Case $\tilde s_{00} = \tilde s_{10} = 0$.} In that case, we have 
$s_{01} = s_{11} = 0$, which is of course
dominated by $\tilde s_{01}, \tilde s_{11}$.

\noindent 2) \textbf{Case $\tilde s_{00} = \tilde s_{10} = 1$.} 
We use corollaries~\ref{cor:survival} and~\ref{cor:boosting}, which we apply to stars $1$ and $3$
independently. Consider star $1$. We first go through a ``boosting'' phase
where we let $N^{(1)}$ drift from $d/4-2$ to $d/3-2$. Then we compute
the probability that $N^{(1)}$ stays above $d/4-2$ for the remaining time. 

\noindent\textbf{Phase 1.} For the boosting phase, we apply Corollary~\ref{cor:boosting}. 
The probability of remaining
below $d/3 - 2$ is at most $2^{-d/21}$.

\noindent\textbf{Phase 2.} The time remaining after boosting is of course at most $M$. 
In time $M$, there is a Poisson
number of steps, say $Q'$, with mean $dM$ (including the steps
where nothing happens). 
From the moment generating
function of the Poisson distribution (see e.g.~\cite{Du96}), we have the following
\begin{equation*}
\prob[Q' \geq d^2M^2]
= \prob[e^{Q'} \geq e^{d^2M^2}]
\leq \frac{\expec[e^{Q'}]}{e^{d^2M^2}}
\leq \frac{e^{dM(e - 1)}}{e^{d^2M^2}}
\leq 2^{-d^2M^2/2}.
\end{equation*}
Assuming $d/3 - 2$
was reached and that there remains at most $d^2 M^2$ discrete steps,
we get that there are at most $d^2M^2$ crossings of the
interval $[d/4 - 3, d/3 - 2]$ by the process $N^{(1)}$. 
By Corollary~\ref{cor:survival},
every time $N^{(1)} = d/3 - 3$, there is a probability of
at least $1 - 2^{-d/12}$ of coming back to $d/3 - 2$ before
hitting $d/4 - 3$. The probability that any of
$d^2 M^2$ attempts at crossing $[d/4 - 3, d/3 - 2]$ succeeds is at most
at most $d^2M^2 2^{-d/12}$ which implies
\begin{equation*}
\prob[\tilde s_{10} = 0] 
\leq d^2 M^2 2^{-d/12} + 2^{-d^2M^2/2} + 2^{-d/21}
\leq 2^{-d/22},
\end{equation*}
for $d$ large enough. Stochastic domination
of the oriented percolation follows directly.

\noindent 3) \textbf{Case $\tilde s_{00} = 1$, $\tilde s_{10} = 0$.} (The
symmetric case is analyzed similarly.) We divide the time window in two phases. 
For the first phase, we compute the probability that
defection ``spreads'' from star $1$ to star $3$. For the second phase, 
we compute the probability that stars $1$ and $3$ remain in or reach
state $1$ respectively. 

\noindent \textbf{Phase 1.} It is easy to see that, in any initial
configuration satisfying $\tilde s_{00} = 1$, $\tilde s_{10} = 0$, six steps (or less)
suffice to reach a configuration with $N^{(3)} \geq 3$. The probability that the
first six steps taken by IPD satisfy this property---call that event $B$---is at least $1/d^6$.
Let $Q''$ be the number of steps until time $M/2$. Then, 
\begin{equation*}
\prob[Q'' \leq 5] \leq 2^{-M/4},
\end{equation*}
by a calculation similar to that in Lemma~\ref{lem:star}.

\noindent\textbf{Phase 2.} We condition on $\{Q'' \geq 6\}$. 
Consider first star $1$. Whether or not $B$ is realized, at the
beginning of Phase 2, we have $N^{(1)} \geq d/4 - 8$. We are back
in the situation of Case 2), except that the time left is only
at least $M/2$. By the same calculation, we obtain that
the probability that $\tilde s_{10}$ is $0$ is at most
$2^{-d/22}$ for $d$ large enough. Consider now star $3$. Let
$Q'''$ be the number of discrete steps left on star $3$. The time
remaining is at least $M/2$. It follows from Corollary~\ref{cor:birth}
that $N^{(3)}$ reaches $d/3 - 2$ before the end of the time
window with probability at least $1/3$ for $d$ large enough.
Once $d/3 - 2$ is reached, we are back to Phase 2 of Case 2).
It follows that on $\{Q'' \geq 6\}$ the probability that
$\tilde s_{11} = 1$ is at least $d^{-6}/4$. Note that
on $\{Q'' \geq 6\}$, the bounds on star $1$ and $3$ are independent.
It is then easy to check that stochastic domination of the
oriented percolation holds.

\noindent$\blacksquare$

\begin{figure}
\begin{center}
\input{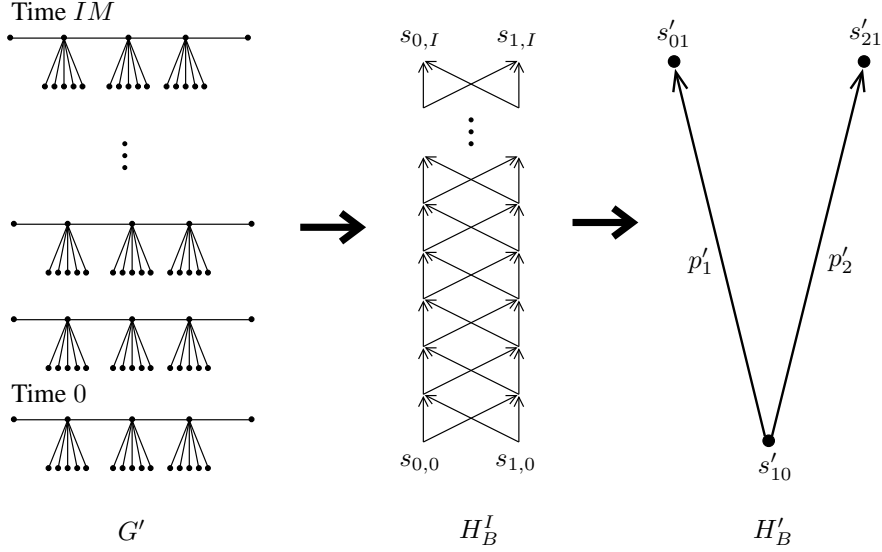}\caption{Further reduction.}
\label{fig:perc2}
\end{center}
\end{figure}

We further simplify the chain by stacking up the construction
in the previous lemma and projecting once more to a smaller space. 
For this, we consider a different percolation model on $\intgr^2$. 
See Figure~\ref{fig:perc2}. Let
$H'_B$ be the directed graph made of three nodes 
$v'_{10} = (1,0), v'_{01} = (0,1), v'_{21} = (2,1)$ with
two edges $e'_1 = (v'_{10}, v'_{01})$, $e'_2 = (v'_{10}, v'_{21})$. 
The edges are open with probability $p'_{1}$, $p'_2$ respectively.
The nodes have state $s'_{10}, s'_{01}, s'_{21}$ respectively
with value in $\{0,1\}$. The percolation works as before with
state $1$ ``traveling'' along open edges. 

Consider again IPD on an arbitrary triplet of stars $G'$ of $G$.
Redefine the vector $\tilde\bfs$ by taking instead
\begin{equation*}
\tilde\bfs = (\tilde s_{00}, \tilde s_{01}, \tilde s_{10}, \tilde s_{11})
= \left(\sigma_d[N_0^{(1)}], \sigma_d[N_{IM}^{(1)}], \sigma_d[N_0^{(3)}], \sigma_d[N_{IM}^{(3)}]\right),
\end{equation*}
for some $I \in \nintgr$ and $M$ as in Lemma~\ref{lem:triplet}. 
We use the following notation: $a\lor b = \max\{a,b\}$.
\begin{lemma}[Towers]\label{lem:tower}
Consider the IPD chain $\{\action_t\}_{t\geq 0}$ on $G=\bush_d^n$ with $d > 15$.
Let $G'$ be an arbitrary triplet of adjacent stars in $G$.
Let $M = d^6$, $I = d^{100}$, and $p'_1 = p'_2 = 1 - 2^{-d/100}$.
Then, for any initial configuration and $s'_{10}$
such that $s'_{10} = \tilde s_{00} \lor \tilde s_{10}$, we have that
$(\tilde s_{01}, \tilde s_{11})$ stochastically
dominates $(s'_{01}, s'_{21})$ for $d$ (constant) large enough. Moreover, the domination holds
for any number of (edge-)nonintersecting triplets simultaneously independently
from each other.
\end{lemma}
\noindent {\bf Proof:} 
The argument ignores any event outside $G'$. 
The proof works by stacking up $I$ copies of $H_B$ and applying
Lemma~\ref{lem:triplet}.
Consider again $\intgr^2$.
We define a \emph{I-tower}, denoted $H_B^I$, to be the graph on nodes
$\{v_{0,i} = (0,i), v_{1,i} = (1,i)\}_{i=0}^I$ where each set of four nodes
of the form $\{v_{0,i}, v_{1,i}, v_{0,i+1}, v_{1,i+1}\}$ induces
a copy of $H_B$ with the same values of $p_0, p_1, p_{10}, p_{01}$
as in Lemma~\ref{lem:triplet}. The node states are denoted
$\{s_{0,i} = (0,i), s_{1,i} = (1,i)\}_{i=0}^I$. By applying repeatedly Lemma~\ref{lem:triplet},
we get that, if $(\tilde s_{00}, \tilde s_{10}) = (s_{0,0}, s_{1,0})$, then
$(\tilde s_{01}, \tilde s_{11})$ stochastically dominates $(s_{0,I}, s_{1,I})$,
so it suffices to show that the latter dominates $(s'_{10}, s'_{21})$.

The case $s'_{10} = 0$ is trivial. So assume $s'_{10} = 1$. Then, the subcase
$\tilde s_{00} \land \tilde s_{10} = 1$ dominates
the subcase $\tilde s_{00} \land \tilde s_{10} = 0$ so it suffices
to consider the latter. Without loss of generality, let 
$\tilde s_{00} = 1$ and $\tilde s_{10} = 0$. The probability
that at least one upwards edge in $H_B^I$ is closed is at most
\begin{equation*}
2I\left(2^{-d/30}\right) \leq 2^{-d/31},
\end{equation*} 
for $d$ large enough.
The probability that no up-right
edge is open is at most
\begin{equation*}
\left(1 - \frac{1}{d^{10}}\right)^{I} \leq 2^{-d/31}, 
\end{equation*}
for $d$ large enough. Therefore,
\begin{equation*}
\prob[s_{0,I} = s_{0,I} = 1] \geq 1 - 2^{-d/32},
\end{equation*}
for $d$ large enough. But note that
\begin{eqnarray*}
\prob[s'_{01} = s'_{21} = 0] 
= (2^{-d/100})^2
= 2^{-d/50}
\geq 2^{-d/32}.
\end{eqnarray*}
So we have domination.

\noindent$\blacksquare$

\subsection{Oriented Percolation}

We conclude the proof of Theorem~\ref{thm:tree} by showing that
the IPD chain at intervals of time $IM$ dominates
a standard percolation model and that in turn the latter model
percolates at an exponential distance from its bottom nodes.

\begin{figure}
\begin{center}
\input{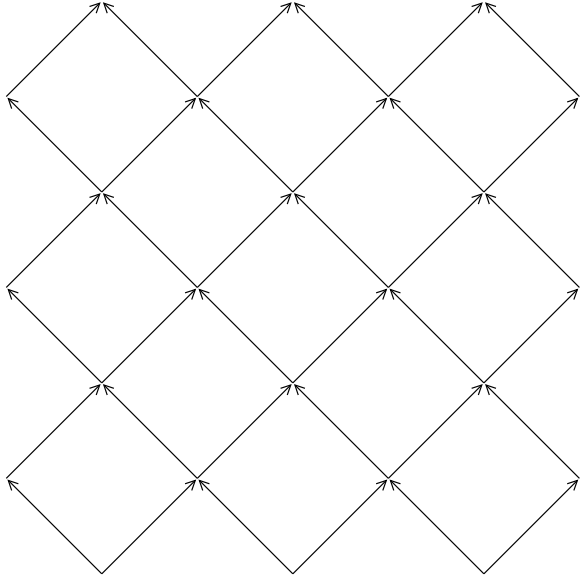}\caption{A section of the oriented percolation lattice.}
\label{fig:perc3}
\end{center}
\end{figure}

\begin{figure}
\begin{center}
\input{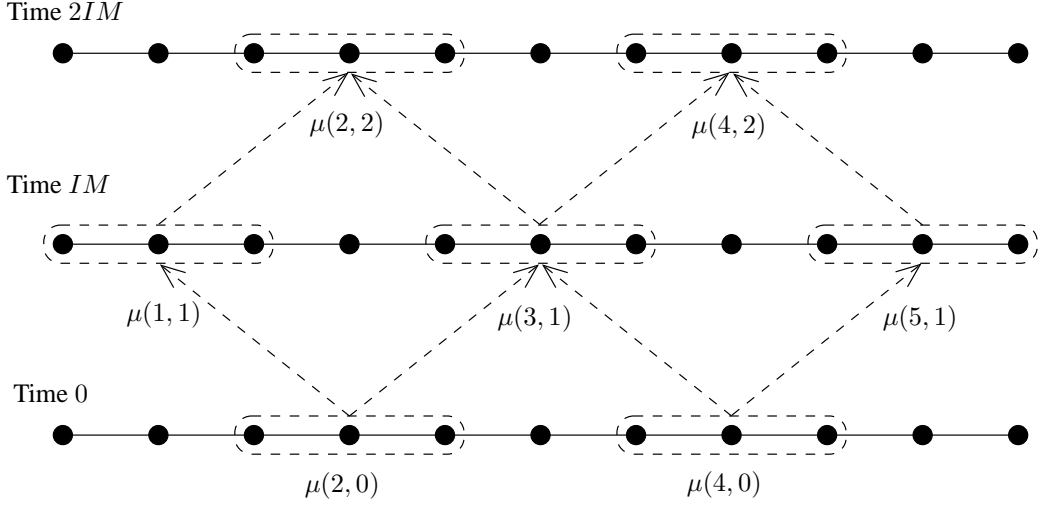}\caption{Graph $\bush_d^{11}$ (stars not shown). Here
$n' = 5$. Circled
triplets correspond to nodes of the percolation lattice of Figure~\ref{fig:perc3}.}
\label{fig:period}
\end{center}
\end{figure}

For convenience, assume $n$ is of the form
\begin{equation*}
n = 2n' + 1,
\end{equation*}
for some positive integer $n'$. (The reason for this
choice will be clear below. See also Figure~\ref{fig:period}.)
Consider the following sublattice of $\intgr^2$,
\begin{equation*}
\pcal = \{(i,j)\in \intgr^2\ :\ 1 \leq i\leq n',\ 0\leq j\leq T,\ i + j\mathrm{\ is\ even}\},
\end{equation*}
where $T$ is a positive integer that will be fixed below. 
Consider the directed graph $G_\pcal = (V_\pcal, E_\pcal)$ with node set
$V_\pcal = \{v_{i,j}\}_{(i,j)\in\pcal}$ and edge set
\begin{equation*}
E_\pcal = \{(v_{i,j}, v_{i+1,j+1}), (v_{i,j}, v_{i-1,j+1})\}_{(i,j)\in\pcal}.
\end{equation*}
See Figure~\ref{fig:perc3} for an illustration. Each edge has probability
$p'$ of being open where $p'$ is set below.
We consider the percolation process on $G_\pcal$ and denote the states
$\bfs'_\pcal = \{s'_{i,j}\}_{(i,j)\in\pcal}$.

Let $\{\action_t\}_{t\geq 0}$ be the IPD chain on
$\bush_d^{n}$ and denote $N^{(i)}_t$ the number
of $\defec$'s on star $i$ at time $t$, excluding
the external nodes.
We consider the following projection of $\{\action_t\}_{t\geq 0}$. 
Let
\begin{equation*}
\mu(i,j) = 4(i-1) + \ind_{\{j\mathrm{\ is\ odd}\}},
\end{equation*}
and let
$\tilde \bfs = \{\tilde s_{i,j}\}_{(i,j)\in\pcal}$ where
\begin{equation*}
\tilde s_{i,j} = \sigma_d\left[N_{jIM}^{(\mu(i,j) - 1)}\right]\lor\sigma_d\left[N_{jIM}^{(\mu(i,j) + 1)}\right],
\end{equation*}
where $I$ and $M$ are as in Lemma~\ref{lem:tower}.
See Figure~\ref{fig:period}. 
We show first
that $\tilde \bfs$ dominates $\bfs'$.
\begin{lemma}[Domination of Oriented Percolation]\label{lem:perc}
Consider the IPD chain $\{\action_t\}_{t\geq 0}$ on $G=\bush_d^n$ with $d > 15$.
Let $M = d^6$, $I = d^{100}$, and $p' = 1 - 2^{-d/100}$.
Let $\action_0 = \alldefec$ (the all-$\defec$ state) and 
let $s'_{i,0} = 1$ for all even $i$'s.
Then, we have that
$\tilde \bfs$ stochastically
dominates $\bfs'$ for $d$ (constant) large enough.
\end{lemma}
\noindent {\bf Proof:} 
This actually follows immediately from
Lemma~\ref{lem:tower}.

\noindent$\blacksquare$

\begin{figure}
\begin{center}
\input{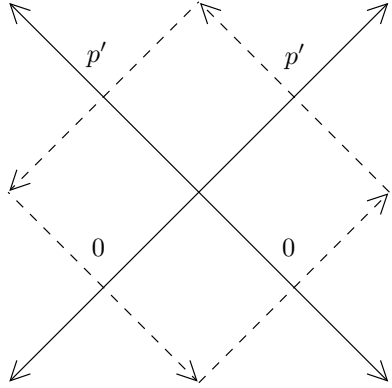}\caption{Dual edges (dashed).}
\label{fig:dual}
\end{center}
\end{figure}

Finally, the next lemma concludes the proof of Theorem~\ref{thm:tree}.
\begin{lemma}[Crossing]
Let $\bfs'$ be defined as above with $p' = 1 - 2^{-d/100}$ and 
let $s'_{i,0} = 1$ for all even $i$'s. Let $T = 2^{(d/2000)n}$.
Assume that $n = 2n' + 1$ and that $T$ is even. Then
\begin{equation*}
\prob[s'_{i,T} = 0,\ \forall i\in\{2,4,\ldots,n'-1\}] \leq
2^{-(d/1000)n},
\end{equation*}
for $d$ (constant) large enough.
\end{lemma}
\noindent {\bf Proof:} 
We use a standard duality argument. For more details, see~\cite{Du84}.
First we modify the percolation lattice $G_\pcal$, which
we now call the primal lattice and still denote $G_\pcal$. To each edge, we add another edge,
reversed, with associated probability of being open $0$. 
We now define the dual lattice. Let
\begin{equation*}
\dcal = \{(i,j)\in \intgr^2\ :\ 1 \leq i\leq n',\ 0\leq j\leq T,\ i + j\mathrm{\ is\ odd}\}.
\end{equation*} 
Consider the directed graph $G_\dcal = (V_\dcal, E_\dcal)$ with node set
$V_\dcal = \{v_{i,j}\}_{(i,j)\in\dcal}$ and edge set
\begin{equation*}
E_\dcal = \{(v_{i,j}, v_{i+1,j+1}), (v_{i,j}, v_{i-1,j+1}),
(v_{i,j}, v_{i-1,j-1}), (v_{i,j}, v_{i+1,j-1})\}_{(i,j)\in\dcal}.
\end{equation*}
Superimpose $G_\pcal$ on top of $G_\dcal$
and notice that to each edge of $G_\dcal$ corresponds
an edge of $G_\pcal$ which is rotated $90^{\mathrm{o}}$ clockwise.
See Figure~\ref{fig:dual}.  
We couple the two lattices
so that an edge in $G_\dcal$ is closed if and only if
the corresponding edge in $G_\pcal$ is open. It is not
hard to see that there is an open path from level $0$ to 
level $T$ in $G_\pcal$ if and only if there is 
no open path from the right boundary to the left
boundary in $G_\dcal$. So it remains to 
compute an upper bound on the latter probability.
Fix any two boundary nodes in $G_{\dcal}$,
say $v_l = v_{1,\eta}$ and $v_r = v_{n',\zeta}$
for some $\eta, \zeta$. The number
of paths of length $L$ between $v_r$ and
$v_l$ is at most $3^L$. Each such path makes
$n' - 1$ more moves to the left than to the right.
In particular, the number of moves to the left
is at least $L/2$. Moreover, each edge going
to the left has a probability $1 - p'$ of 
being open. So the probability that there
is a path between $v_r$ and $v_l$
(which we denote $v_r \to v_l$) is at most
\begin{eqnarray*}
\prob[v_r \to v_l]
\leq 
\sum_{L=n'-1}^{+\infty} 3^L (1 - p')^{L/2}
\leq \frac{(3 \, 2^{-d/200})^{\frac{n-1}{2}}}{1 - 3\, 2^{-d/200}},
\end{eqnarray*}
for $d$ large enough. There are at most $T^2$ pairs of boundary
nodes so by the union bound
\begin{eqnarray*}
\prob[s'_{i,T} = 0,\ \forall i\in\{2,4,\ldots,n'-1\}] 
\leq T^2 \frac{(3 \, 2^{-d/200})^{\frac{n-1}{2}}}{1 - 3\, 2^{-d/200}}
\leq 2^{-(d/1000)n},
\end{eqnarray*}
for $d$ large enough.

\noindent$\blacksquare$

\section{Win-Stay Lose-Shift on Graphs with Large Expansion}\label{sec:expand}

For this section, we consider the discrete-time version of the chain.
That is, at every time step, we pick one edge uniformly at random and
update the actions at the endpoints of that edge. Equivalently,
we look at the discrete-time chain embedded in $\{\action_t\}_{t\geq 0}$
by stopping the chain every time a clock rings. Also, since
we are looking for a lower bound on $T_\allcoop$, we can speed
up the chain by picking only those edges with at least 
one $\defec$ endpoint. Denote the discrete-time sped-up
chain $\{\baction_k\}_{k \in \nintgr}$. 

The proof of Theorem~\ref{thm:expand} is based on the following
geometric observation. Let $U_k$ be the set of nodes defecting at time $k$ and
denote $N_k = |U_k|$. At the next update, $N_k$ goes down by 2 if we pick an edge
``inside'' $U_k$ and it goes up by 1 if we pick an edge on the ``boundary'' of $U_k$.
Therefore, if the boundary of $U_k$ is more than twice as big as the inside
of $U_k$, on average the chain moves away from the fixed point $\allcoop$. 

\noindent \textbf{Proof of Theorem~\ref{thm:expand}:}
Let $U \subseteq V$ with $\alpha(n) \leq |U| \leq \beta(n)$.
Note first that $\rho_{\alpha,\beta}(G) > 1/2 + \eps$ implies
\begin{equation*}
|E(U, U^c)| \geq \left(\frac{1}{2} + \eps\right) \vol(U).
\end{equation*}
Let $\eps' > 0$ such that $2 - \eps' = (1/2 + \eps)^{-1}$. Then
\begin{equation*}
|E(U, U^c)| + 2|E(U,U)| = \vol(U)
\leq (2 - \eps')|E(U,U^c)|,
\end{equation*}
which implies 
\begin{equation*}
2 |E(U,U)| \leq (1 - \eps')|E(U,U^c)|.
\end{equation*}
Therefore there is an $\eps'' > 0$ such that
if $\alpha(n) \leq N_k \leq \beta(n)$, then
\begin{displaymath}
N_{k+1} 
= \left\{\begin{array}{ll}
N_{k} + 1, & \mathrm{with\ probability\ at\ least\ }\frac{2}{3} + \eps'',\\
N_{k} - 2, & \mathrm{with\ probability\ at\ most\ }\frac{1}{3} - \eps''.
\end{array}
\right.
\end{displaymath}

Let 
\begin{equation*}
a = \left[\frac{1}{2}\left(\frac{2}{3} + \eps''\right)
\left(\frac{1}{3} - \eps''\right)^{-1}\right]^{1/3} > 1.
\end{equation*}
It is easy to check that
\begin{equation*}
\left(\frac{2}{3} + \eps''\right) a^{-1} +
\left(\frac{1}{3} - \eps''\right) a^{2} < 1.
\end{equation*}
Therefore,
\begin{equation*}
W(N_k) = a^{n - N_k},
\end{equation*}
is a bounded nonnegative supermartingale on $\{\alpha(n) \leq N_k \leq \beta(n)\}$.
Using the optional sampling theorem as in~\cite{DG+02},
it follows that the probability of $N_k$ crossing the interval $[\alpha(n), \beta(n)]$
is less than $a^{-(\beta(n) - \alpha(n))}$ for $n$ large enough. The theorem immediately follows.

\noindent $\blacksquare$

\section{Concluding Remarks}

The work presented here leads naturally to the following questions:
\begin{enumerate}
\item Is there a $d$ (constant) such that for all
$n$ large enough and for all
trees of minimum degree $d$ with $n$ nodes,
the emergence of cooperation is exponentially slow?
\item What is a good criterion for fast emergence
of cooperation in this setup? Is the line
and its---appropriately defined---variants the only graphs on which
the convergence to all-cooperation is fast?
\end{enumerate}

\section*{Acknowledgments}

The first author acknowledges the support of a 
Miller Fellowship in Statistics and Computer Science, U.C. Berkeley, 
a Sloan Fellowship in Mathematics and NSF grants DMS-0504245 and
DMS-0528488. 
The second author is supported by CIPRES (NSF ITR grant \# NSF EF 03-31494), FQRNT, NSERC 
and a Lo\`eve Fellowship. The second author also thanks Martin Nowak and the
Program for Evolutionary Dynamics at Harvard where part of this work was done.


\clearpage

\begin{thebibliography}{99}

\bibitem[Ax84]{Ax84}
Axelrod, R. (1984).
\textit{The Evolution of Cooperation}.
Basic Books.

\bibitem[BK$^+$05]{BKMP05}
N.~Berger, C.~Kenyon, E.~Mossel, and Y.~Peres.
\newblock Glauber dynamics on trees and hyperbolic graphs.
\newblock {\em Probab. Theory Related Fields}, 131(3):311--340, 2005.
\newblock Extended abstract by Kenyon, Mossel and Peres appeared in proceedings
  of 42nd IEEE Symposium on Foundations of Computer Science (FOCS) 2001,
  568--578.

\bibitem[Du84]{Du84}
Durrett, R. (1984).
Oriented Percolation in Two Dimensions.
\textit{The Annals of Probability}
\textbf{12} 999--1040.

\bibitem[Du96]{Du96}
Durrett, R. (1996).
\textit{Probability: Theory and Examples}.
Duxbury.

\bibitem[DG$^+$02]{DG+02}
Dyer, M., Goldberg, L.A., Greenhill, C.,
Istrate, G., and Jerrum, M. (2002).
Convergence of the Iterated Prisoner's
Dilemma Game. \textit{Combinatorics, 
Probability, and Computing} 
\textbf{11} 135--147.

\bibitem[FKS89]{FKS89}
Friedman, J., Kahn, J., and Szemeredi, E. (1989).
On the second eigenvalue in Random Regular Graphs.
In: \textit{Proceedings of ACM STOC}.

\bibitem[FL98]{FL98}
Fudenberg, D. and Levine, D.K. (1998).
\textit{The Theory of Learning in Games}.
MIT Press.

\bibitem[Ka95]{Ka95}
Kahale, N. (1995). 
Eigenvalues and Expansion of Regular Graphs.
\textit{Journal of the ACM},
\textbf{42} 1091--1106. 

\bibitem[Ki95]{Ki95}
Kittock, J.E. (1995). Emergent conventions and the structure
of multi-agent systems. In: 
\textit{1993 Lectures in Complex Systems: Proceedings
of the 1993 Complex Systems Summer School},
Nadel, L. and Stein, D, eds. Vol. VI of
\textit{Santa Fe Institute Studies in the Sciences
of Complexity Lecture}, Santa Fe Institute,
Addison-Wesley.

\bibitem[Li85]{Li85}
Liggett, T.M. (1985).
\textit{Interacting Particle Systems}.
Springer.

\bibitem[Li99]{Li99}
Liggett, T.M. (1999).
\textit{Stochastic Interacting Systems: Contact, Voter and Exclusion Processes}.
Springer.


\bibitem[NS93]{NS93}
Nowak, M. and Sigmund, K. (1993). A strategy of win-stay, lose-shift
that outperforms tit-for-tat in the Prisoner's Dilemma game. \textit{Nature}
\textbf{364} 56--58.

\bibitem[ST93]{ST93}
Shoham, Y. and Tennenholtz, M. (1993). Co-learning and
the evolution of social activity. \textit{Mimeo}.

\bibitem[ST97]{ST97}
Shoham, Y. and Tennenholtz, M. (1997). On the emergence of social
conventions: Modelling, analysis and simulations. \textit{Artificial Intelligence}
\textbf{94} 139--166.

\end{thebibliography}
\end{document}